\documentclass[letterpaper, 10pt, conference]{ieeeconf}      
\IEEEoverridecommandlockouts                              
                                                          
\let\labelindent\relax
\usepackage{graphicx}   
\usepackage{epsfig} 
\usepackage[T1]{fontenc} 
\usepackage{amssymb} 
\usepackage{amsmath} 
\usepackage{algorithm}
\usepackage[noend]{algpseudocode}
\usepackage{graphicx}
\usepackage{lipsum}%
\usepackage{xcolor}
\usepackage{hyperref}
\usepackage{float}
\usepackage{enumitem}
\usepackage{verbatim}
\usepackage{comment}
\usepackage{booktabs}
\usepackage{svg}
\usepackage[backend=bibtex,style=ieee,sorting=none]{biblatex}
\usepackage{subcaption}
\usepackage{mwe}

\usepackage{color}

\newcommand{\rom}[1]{\uppercase\expandafter{\romannumeral #1\relax}}

\bibliography{SourcesAll.bib}

\title{\LARGE \bf Informed sampling-based trajectory planner for automated driving in dynamic urban environments }

\author{Robin Smit$^{1}$, Chris van der Ploeg$^{1,2}$, Arjan Teerhuis$^1$, Emilia Silvas$^{1,3}$
\thanks{$^{1}$Netherlands Organisation for Applied Scientific Research, Integrated Vehicle Safety Group, 5700 AT Helmond, The Netherlands.}
\thanks{$^{2}$Eindhoven University of Technology, Dynamics and Control Group, Mechanical Engineering Dept., P.O. Box 513, 5600 MB, Eindhoven, The Netherlands.}
\thanks{$^{3}$Eindhoven University of Technology, Control Systems Technology Group, Mechanical Engineering Dept., P.O. Box 513, 5600 MB, Eindhoven, The Netherlands.}
}

\date{April 2021}
\begin{document}
\maketitle
\thispagestyle{empty}
\pagestyle{empty}
\begin{abstract}
The urban environment is amongst the most difficult domains for autonomous vehicles. The vehicle must be able to plan a safe route on challenging road layouts, in the presence of various dynamic traffic participants such as vehicles, cyclists and pedestrians and in various environmental conditions. The challenge remains to have motion planners that are computationally fast and that account for future movements of other road users proactively. This paper describes an computationally efficient sampling-based trajectory planner for safe and comfortable driving in urban environments. The planner improves the Stable-Sparse-RRT algorithm by adding initial exploration branches to the search tree based on road layout information and reiterating the previous solution. Furthermore, the trajectory planner accounts for the predicted motion of other traffic participants to allow for safe driving in urban traffic. Simulation studies show that the planner is capable of planning collision-free, comfortable trajectories in several urban traffic scenarios. Adding the domain-knowledge-based exploration branches increases the efficiency of exploration of highly interesting areas, thereby increasing the overall planning performance. 
\end{abstract}
\section{Introduction}\label{sec:introduction}
Since dense traffic and cluttered environments bring uncertainties and interaction challenges with other road users, one of the most challenging tasks in automated driving is safe and comfortable urban driving. Amongst others, trajectory planning is a far from trivial task due to the varying road layouts, different (types of) road users, and highly dynamic scenarios. The task of planning a trajectory for the automated vehicle (AV) to follow has gained increasing attention and focus over the last years \cite{Paden2016}\cite{Claussmann2019}.

To compute fast (near-)optimal trajectories in terms of safety and comfort, several types of algorithms exist, falling in one of these categories: \textit{optimization-based}, \textit{path primitive-based} and \textit{sampling-based}. The first approach formulates the environment, such as the road layout and other traffic participants, as artificial potential fields or hard constraints and tries to solve the optimization problem analytically. For instance, \cite{Siebenrock2020} tries to find an optimal trajectory by formulating a Model Predictive Control (MPC) problem. However, usually such methods are limited in flexibility of environment representations due to the analytical formulation, limiting their applicability in real-world autonomous driving systems. The second approach generates a set of path primitives and calculates an optimal velocity and/or acceleration profile per path primitive, after which the most optimal trajectory in terms of safety and comfort is selected. In \cite{Artunedo2019}, a real-time motion planner is developed based on sampling path primitives as $5^{th}$-order Bezier curves and calculating suitable velocity profiles based on the motion prediction of dynamic objects. However, the fixed mathematical description of the path primitives can limit the flexibility of the solution, especially in complex road layouts and highly cluttered environments. Moreover, the decoupling of the velocity profile from the path primitives causes problems for online decision making such as overtaking a slower vehicle versus cruising behind it. 

To overcome the limited flexibility of the aforementioned approaches, the trajectory planning problem can be solved by using a \textit{sampling-based} approach, which has shown great potential in several robotics applications \cite{Arab2016}\cite{Littlefield2018}. In \cite{Li2014}, an algorithm called \textit{Stable-Sparse-RRT (SST)} is introduced where piecewise-constant control samples are generated following a stochastic distribution and propagated through a kinematic model. Sampling and propagating control inputs guarantee that the solution adheres to the vehicles non-holonomic kinematic constraints. Furthermore, by applying forward propagation of the kinematic system in time, a time space can be mapped on the system states. A time-bounded collection of vehicle states will be henceforth referred to as a trajectory. Mapping calculated vehicle states in time enables the planner to take time-bounded motion prediction of other road users into account. However, due to the inherently exploring nature of such sampling based algorithms, convergence to high quality trajectories is slow and (sub-)optimality is often not guaranteed \cite{Zheng2021}. To overcome these limitations and enable SST solutions to be applied for AV trajectory planning, initial state space exploration of high potential regions such as the lane center or around a previously found solution is enforced, resulting in higher quality solutions while maintaining the same planning update rate. 

The main contributions of this work are as follows: 
\begin{enumerate}
    \item We apply the generic algorithm introduced in \cite{Li2014} in a real-time automated driving trajectory generation application 
    \item We present a novel extension to the trajectory planner, improving its efficiency, by providing an initial exploration method of high potential regions such as the lane center and previous solutions.
\end{enumerate}
The proposed trajectory planner is implemented in a realistic simulation environment and validated based on several urban-driving scenarios, showing its capability of providing safe and comfortable trajectories in the case of complex road layouts and other static and dynamic traffic participants. 

This article is outlined as follows. In Section \ref{sec:problemstatement}, the problem statement and preliminaries are defined. Section \ref{sec:method} focuses on the method used for solving the trajectory planning problem. This method is proven in a simulation study in Section \ref{sec:validation} and future work is presented in Section \ref{sec:conclusion}.

\section{Preliminaries}\label{sec:problemstatement}
To improve on-road safety in extended operational design domains, novel motion planning algorithms are developed in the EU Horizon Safe-up Project. Based on prediction models of other traffic participants' behaviour and advanced risk estimation, motion planners and vehicle control strategies need to show improved safety and capabilities in various use cases, such as bad weather, faulty sensors or dynamic urban environments including pedestrians or cyclists. 

The focus of this work is to find safe, comfortable and feasible vehicle trajectories in an urban environment. Due to the highly dynamic nature of this environment, the trajectory planner should be able to find suitable solutions in real-time with a sufficient update rate. Furthermore, the solution should take into account the kinodynamic constraints of the vehicle, the road layout and the predicted motion of other traffic participants. 
\subsection{Problem statement}
\begin{figure}[h!]
    \centering
    \includegraphics[clip, trim=75mm 60mm 75mm 60mm, width=1\columnwidth]{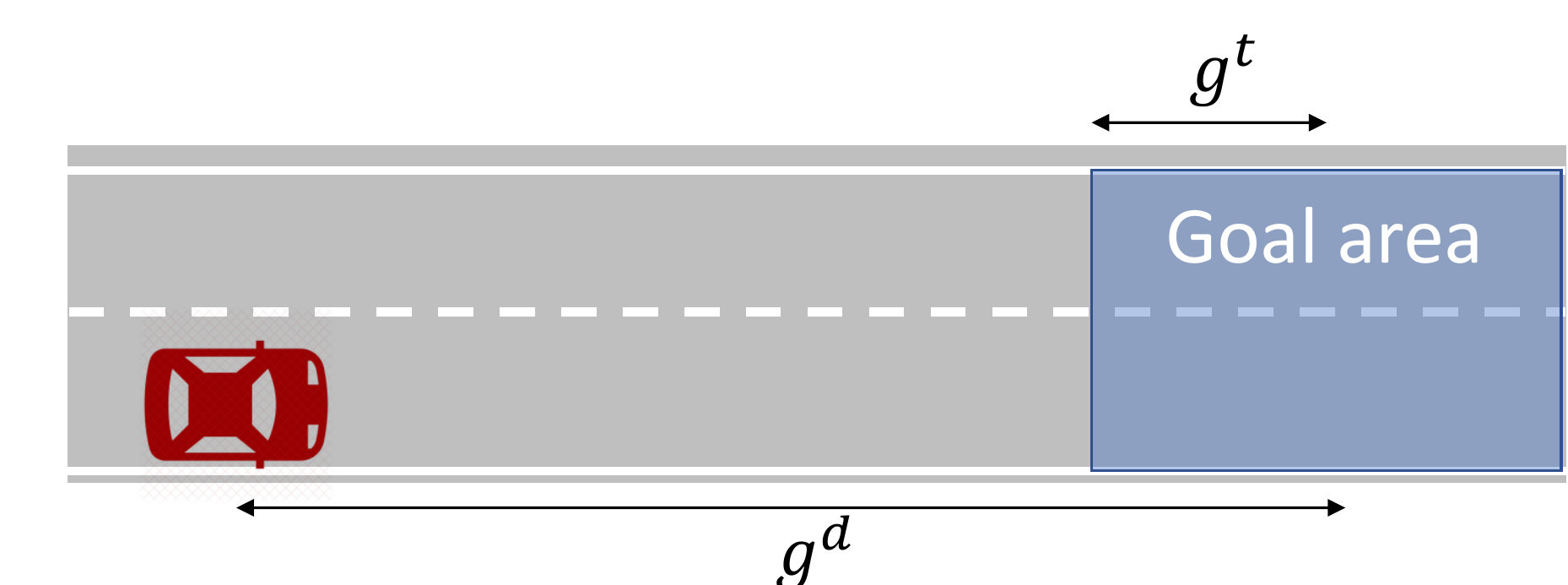}
    \caption{\footnotesize{Ground plane visualization of goal space}}
    \label{fig:goalspace}
\end{figure}

To define the constraints and goals of the planner, first consider the state of the vehicle at time $k$ denoted by $\textbf{x}_k=(x_k,y_k,\theta_k,v_k)$, where $x$ and $y$ represent the vehicle ground-plane position, $\theta$ is the vehicle planar orientation and $v$ is the vehicle longitudinal velocity. To allow for safe and comfortable driving in urban environments, the vehicle needs to proactively account for the predicted motion of other traffic participants. Often these predictions are neglected which can lead to reactive and thus more aggressive vehicle behavior. To account for predicted motions, assume a set of $O$ objects, where the measured state of an object at time $k$ is given by its planar coordinates denoted by $\textbf{o}^o_k=(x^{o}_k, y^{o}_k, \theta^{o}_k), \forall o \in O$. Furthermore, assume the predicted trajectories of these dynamic objects to be available with no noise or inaccuracies. The challenge is then to include these into the trajectory planning problem in a (near) optimal way by balancing progress along the road and minimize driving risks. In the absence of static or dynamic obstacles, the vehicle is expected to follow the lane centers of the road network. To model that, the road infrastructure is described by a set of $L$ lanes, where each lane is described as a collection of lane center points, together with a lane width. A lane is denoted by $l^i = ( r^1, r^2, \hdots, r^R, \omega ), \forall i \in [1 \hdots L]$, where $r^j \forall j \in R$ is a lane center point described by its planar position, $(x^r, y^r)$ and $\omega$ is the lane width. To lead the vehicle along a global route, the trajectory planner receives global goals to plan towards. Here, the goal is represented by a goal space $\mathbb{G} \subset \mathbb{R}^4$, and it is computed such that its position components $(x,y)$ span over all available lanes at a fixed distance $g^d$ along the vehicles position, with longitudinal threshold $g^t$. Furthermore, the orientation component $\theta$ and velocity component $v$ are unbounded, meaning that a specific orientation and velocity are not required for the goal to be reached. The schematic overview of the goal space is depicted in Figure \ref{fig:goalspace}.
Considering the aforementioned requirements for safe and comfortable driving, the trajectory planning problem at time $k$ requires finding a trajectory $ \xi_k  =  \begin{bmatrix}
\textbf{x}_0 \;
\textbf{x}_1 
\hdots
\textbf{x}_f \\
\end{bmatrix}$, which satisfies the kinodynamic constraints of the vehicle, is collision-free given the predicted object states and brings the vehicle to the goal space such that $\textbf{x}_f \in \mathbb{G}$. Besides these constraints, the planner should also minimize the risk of collision for the vehicle as well as other traffic participants by maintaining sufficient distance to others, follow the appropriate lane center for a given road network whenever possible and maximize comfort.
\subsection{State and control space definition}
\begin{figure}
    \centering
    \includegraphics[width=0.8\columnwidth]{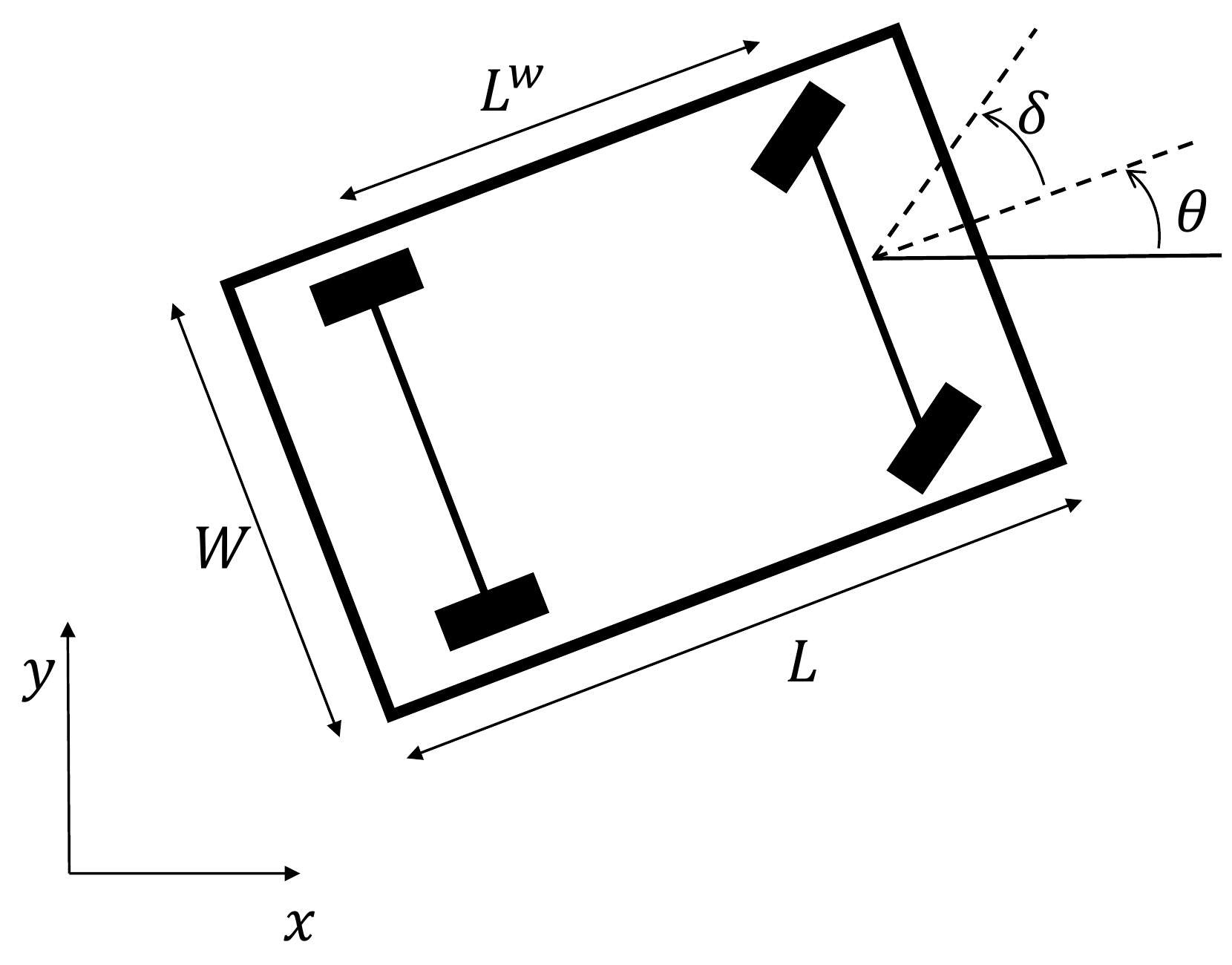}
    \caption{\footnotesize{Schematic representation of the kinematic bicycle model.}}
    \label{fig:kinematicmodel}
\end{figure}

Most sampling-based planners generate search trees by sampling new states in a selected state space. Such state spaces can be $\mathbb{R}^2$, Dubins or Reeds-Shepp \cite{Rathinam2018}. However, planning in such spaces either puts a limit on the movement of the car (e.g., the constant curve radius in Dubins), does not satisfy the differential and/or non-holonomic kinematic constraints of the vehicles motion model (e.g., instantaneous discrete change of vehicle heading) or require a specific steering function for the system. To facilitate the integration of kinematic and differential constraints, while maintaining flexibility on the reachable space of the vehicle, control space sampling poses an alternative. Here, control inputs are sampled and propagated in the system over a specified time. Given a fixed initial state (sample), this results in a new state space sample that satisfies the kinematic constraints. 
As often used in motion planning designs and also in this work, a non-linear kinematic bicycle model as shown in Figure \ref{fig:kinematicmodel} is used to allow for real-time trajectory planning for larger time horizons, i.e., 4-10 seconds \cite{Paden2016}. 
The model is described by its discretized kinematic equations, using Euler discretization, as:
\begin{equation}
\label{eq:vehiclemodel}
\begin{split}
    x_{k+1} =& x_k+t^sv_k\cos{\theta_k},\\
    y_{k+1} =& y_k+t^sv_k\sin{\theta_k},\\
    \theta_{k+1} =& \theta_k+t^s\frac{v_k}{L^w}\tan{\delta_k},\\
    v_{k+1} =& v_k+t^sa_k,
\end{split}
\end{equation}
where $t^s$ is the integration time, $a$ denotes the longitudinal acceleration of the vehicle, $\delta$ denotes the steering angle and $L^w$ denotes the wheel base of the vehicle. 
Due to physical constraints of the vehicle, such as a maximum velocity, acceleration and/or steering angle, the state is bounded by $\textbf{x} \in \mathbb{X} \subset \mathbb{R}^4 $ and the input to the system is denoted by  $\textbf{u} = [ a,  \delta] \in \mathbb{U} \subset \mathbb{R}^2 $. 

\subsection{Stable Sparse RRT (SST)} \label{SST}
\begin{figure}[h!]
    \centering
    \includegraphics[clip, trim=80mm 30mm 80mm 30mm, width=\columnwidth]{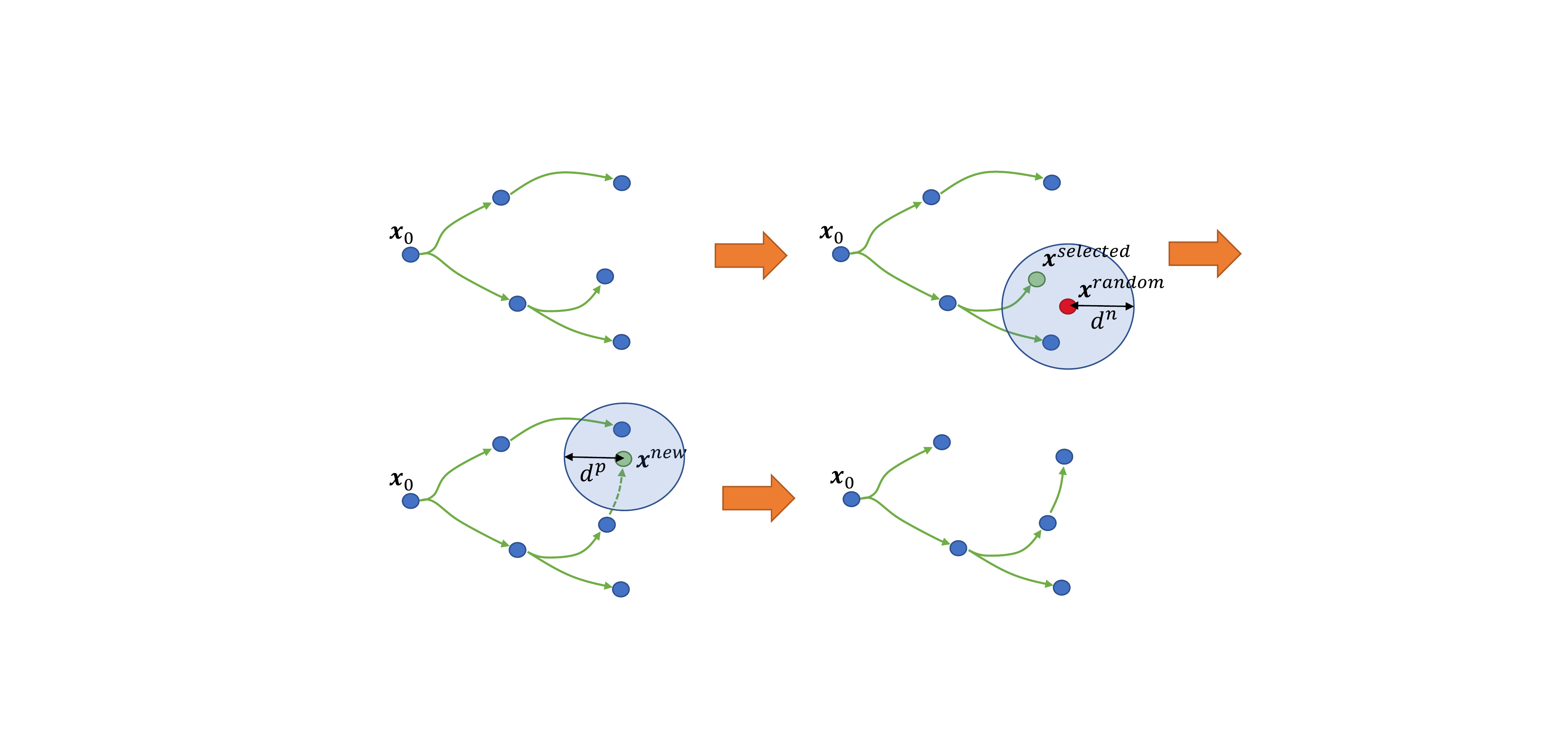}
    \caption{\footnotesize{Simplified overview of the SST \textit{selection-propagation-pruning procedure}}}
    \label{fig:SSTProcedure}
\end{figure}

To plan safe and comfortable trajectories for autonomous urban driving, this work builds upon the anytime-planner introduced in \cite{Li2014}, where a sampling-based planning method called \textit{Stable Sparse RRT (SST)} is introduced. Given an initial state, SST explores the state space by generating a state tree, following an iterative \textit{selection-propagation-pruning} procedure (as shown in Figure \ref{fig:SSTProcedure}), briefly described here and more extensively in \cite{Li2014}. 

\textit{i) Selection}. First, a random state $\textbf{x}$ is sampled in $\mathbb{X}$. Then a set of states in the current tree are selected within distance $d^n$ of this random state. If no states are available within $d^n$, the closest state in the state tree is returned. If more than one state are found in the vicinity of the random state, the state is returned with the lowest motion cost-to-go from the start state. A detailed description of the motion costs applied in this work is given in Section \ref{optimization}. 

\textit{ii) Propagation}. A control input is sampled and the selected state is propagated for propagation time $t^p$ to achieve a new state. Details on the input sampling procedure can be found in Section \ref{inputsampling}. 

\textit{ iii) Pruning}. Finally, the new state is added to the tree only if it passes several validity checks, as described in Section \ref{validitychecking}, and the closest state is either more than threshold $d^p$ away, or the motion cost-to-go is lower than that of the set of states in its vicinity. If the new state is added to the tree, the states in the vicinity are moved to an inactive subset and thus removed from further consideration in the planning iteration. 

\section{Trajectory generation method}\label{sec:method}
To find a fast, safe and comfortable vehicle trajectory, the SST algorithm, developed previously for robotic applications, is applied in this work in an automated driving architecture, and enhanced with specific methods for sampling the control inputs, checking the validity of vehicle states and calculating the motion costs of travelling to successive states. Furthermore, a novel method is proposed expanding the SST algorithm, to force initial exploration of high potential state space regions.
\subsection{Input sampling}\label{inputsampling}
The SST algorithm applies piece-wise constant control inputs to the kinodynamic model to steer the system from a selected state to a new state. Here, a control input sampler is proposed which takes into account the physical limitations of the system in terms of achievable acceleration, deceleration and steering angle by bounding the input samples. The input components are sampled following a zero-mean Gaussian distribution with standard deviation $\sigma^a$ and $\sigma^{\delta}$, for the acceleration and steering angle respectively. If the sample falls outside of the input bounds, $[\underline{a}, \overline{a}]$ and $[\underline{\delta}, \overline{\delta}]$ respectively, a new sample is drawn until a valid input is found. A Gaussian distribution is chosen to promote the probability of samples around the zero-mean, but still allows for steering, braking and accelerating up to the specified limits whenever required. The sampled inputs are applied to the system by forward propagation during time $t^p$. Note that propagation time $t^p$ can be equal to, but also be a multitude of integration step time $t^s$.
\subsection{State validity checking}\label{validitychecking}
To prevent the planner from steering the system to states which are not valid, a state validity checker is proposed in this section. To keep the planner as efficient as possible, non-valid states are not taken into consideration in the search tree generation. According to the trajectory planning problem definition, a state is considered valid if it is within specific state bounds, within the lanes and the full vehicle geometry is not in collision with other traffic participants.  

\textit{i) State bounds:} The validity checker checks for two different types of state bounds. Position state bounds on the $x$ and $y$ state are added to prevent exploration in undesired regions of $\mathbb{X}$, such as positions far beyond the ground plane area of $\mathbb{G}$, thereby stimulating faster convergence to high-quality paths. Furthermore, a bound is added to the velocity state based on the local legal speed limit or vehicle physical limitation. 

\textit{ii) Road layout:} The vehicle should remain within the driveable space, spanned by the lanes, at all times. To check whether a vehicle is within a lane, a penalty grid map $P(x,y)$ is constructed. In this grid map, the cell penalties increase linear along the distance to the closest lane center point, up to a maximum value $\overline{P}$, as follows:
\begin{equation}
    P(x,y) = 
    \begin{cases}
    \frac{2\overline{P} d^{r}}{\omega}, & \text{if } d^{r} < \frac{\omega}{2} \\
    \overline{P}, & \text{otherwise,}
    \end{cases}
\end{equation}
where $d^{r}$ is the Euclidean distance from the cell center to the closest lane center point in $L$.

\begin{figure}[h!]
    \centering
    \includegraphics[clip, trim=10mm 10mm 10mm 20mm, width=0.9\columnwidth]{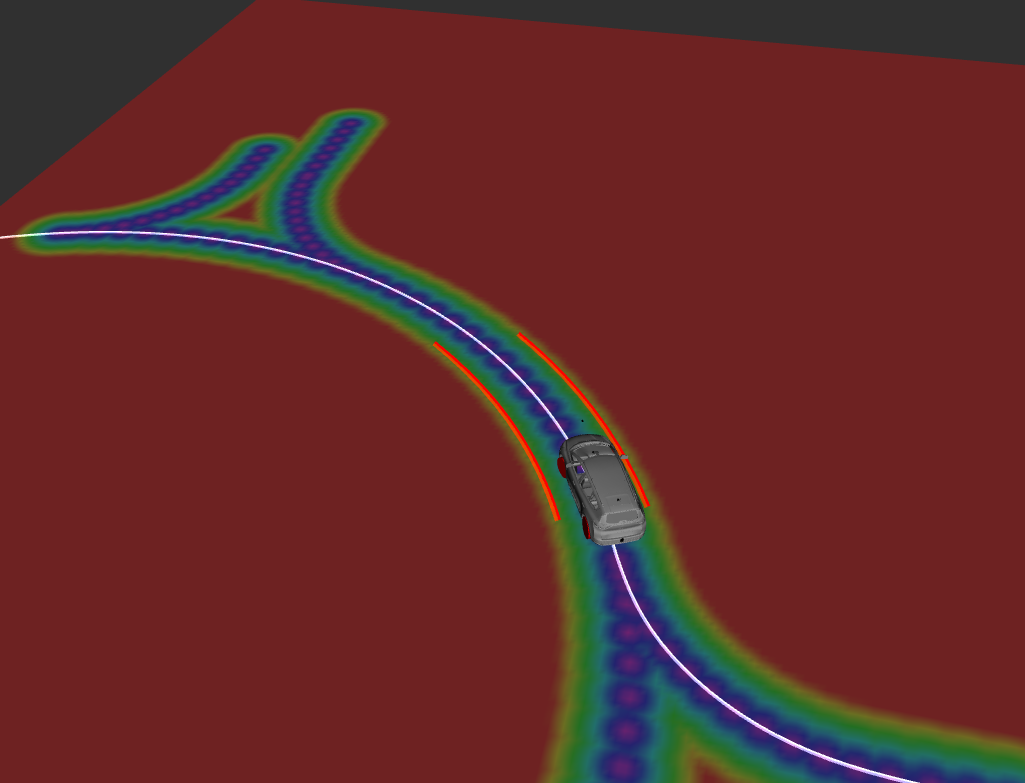}
    \caption{\footnotesize{Grid map based on lane information on a roundabout with red lines representing lane markings.}}
    \label{fig:gridmap}
\end{figure}

The ground plane position part, $(x,y)$, of the vehicle state under consideration is mapped on a cell in a penalty grid map. If the penalty of this cell is equal to threshold $\overline{P}$, the vehicle state is considered invalid. 
An example of the penalty grid map is depicted in Figure \ref{fig:gridmap}, where the lane center is given a low penalty (purple) and the penalty increases towards the lane edges, until it reaches its maximum value $\overline{P}$ (red).

\textit{iii) Object collision:} Any vehicle state where the ego vehicle geometry, denoted by length $L$ and width $W$, overlaps with the bounding box of another traffic participants predicted state at the same time instance is considered as a collision. For a given vehicle pose at time $k$, $(x_k, y_k, \theta_k)$, a polygon is created on the ground plane considering the length and width of the vehicle. Similar polygons are created for all predicted target poses at time $k$, under the assumption that the geometry of any target is considered known. If the vehicle polygon overlaps with at least one of the predicted target polygons, the state is considered invalid. The polygon overlap is checked using the \textit{ray casting} algorithm \cite{Kularathne2018}. 

\subsection{Optimization objectives}\label{optimization}
Due to the \textit{anytime}-nature of the planner, every planning iteration is allowed to run for a fixed amount of time, $t^q$, tuning which is a trade-off between planner update rate and quality of the found trajectory. To determine the most optimal solution, costs are awarded either to states, \textit{state costs}, or the transition from one state to the other, \textit{motion costs}. Resulting motion costs are calculated from state costs of successive states by approximating the integral of the cost using the trapezoidal approximation function. Henceforth, any state cost of state $\textbf{x}$ will be referred to as $c(\textbf{x})$, whereas the motion cost from state $\textbf{x}_n$ to state $\textbf{x}_{n+1}$ will be denoted by $m(\textbf{x}_n,\textbf{x}_{n+1})$. 

The total motion cost of travelling from a starting state $\textbf{x}_0$ to a state $\textbf{x}_N$ is a summation of the successive motion costs: 
\begin{equation}
m(\textbf{x}_0,\textbf{x}_N) = \sum^{N-1}_{n=0} m(\textbf{x}_n,\textbf{x}_{n+1}).
\end{equation}
In this work, each individual motion cost is a weighted summation of four separate optimization objectives: \textit{\textbf{p}ath \textbf{l}ength}, \textit{\textbf{d}esired \textbf{v}elocity}, \textit{\textbf{p}enalty \textbf{g}rid} and \textit{\textbf{t}arget \textbf{c}learance}:
\begin{equation}
    \tilde{m}(\textbf{x}_n, \textbf{x}_{n+1}) = \sum_{i \in W} w^i m^i(\textbf{x}_n, \textbf{x}_{n+1}), W = \{pl, dv, pg, tc\}, 
\end{equation}
where $\tilde{m}(\cdot)$ is the total multi-objective motion cost and $w^i$ is the weighting of the individual optimization objectives, which are detailed in the remainder of this section.
 
\paragraph{Path length}
Where comfort and safety allow it, the vehicle should attempt to reach its goal as efficient as possible in terms of driving distance. Therefore, this optimization objective tries to minimize the geometric length of the trajectory, i.e., the distance along the ground plane path, by adding a cost based on the Euclidean distance between the positions of successive states:
\begin{equation}
    m^{pl}(\textbf{x}_n, \textbf{x}_{n+1}) = \sqrt{(x_{n+1}-x_n)^2 + (y_{n+1}-y_n)^2}.
\end{equation}

\paragraph{Desired velocity} 
In nominal traffic, the vehicle should try to drive as fast as possible while remaining safe, comfortable and obliging to traffic rules. Therefore, an objective is added to push the vehicle velocity towards a desired speed, which is (near) the legal speed limit or at another situation-appropriate speed advice. The state cost of the desired velocity objective is defined as:
\begin{align}
    \label{eq:cds}
    c^{dv}(\textbf{x}_n) = |v_n - v^{d}|,
\end{align}
where the desired velocity is denoted by $v^d$.
\paragraph{Penalty grid}
Besides rendering states invalid, the penalty grid as described in section \ref{validitychecking} also penalizes undesired states, using the cost term: 
\begin{equation}
    c^{pg}(\textbf{x}_n) = P(x_n, y_n).
\end{equation}
This optimization objective makes it attractive for the vehicle to drive on the lane center. When a deviation from the lane center is required, for instance for a static object near the side of the road, this optimization objective ensures that the lateral displacement with respect to the center of the lane is as minimal as safety allows. 

\paragraph{Target clearance}
Due to the dynamic nature of the urban environment, trajectory planning should be able to take into account the predicted motion of other traffic participants. Besides collision avoidance by rendering a state invalid when its position collides with a predicted target position, the vehicle should steer clear of other traffic participants with some margin. Therefore, a repulsive artificial potential field is introduced around the time-dependent target position. The potential field of a perceived object is described using a two-dimensional Gaussian function as follows:
\begin{align}
    c^{tc}(\textbf{x}_k)=&\sum_{n=1}^{N}A^{o}e^{-f(\mathbf{x}_k,\mathbf{o}^o_k)}, \\
    \text{where }f=&\frac{\left(x_k-x^{o}_k\right)^2}{\left(\sigma^{o,x}\right)}+\frac{\left(y_k-y^{o}_k\right)^2}{\left(\sigma^{o,y}\right)},
\end{align}
where $x^o_k$, $y^o_k$ represent respectively the predicted x-position and y-position at time $k$ of object $o$. Furthermore, $A^o$ represents the maximum magnitude of the repulsive field of the object. Finally, $\sigma^{o,x}$ and $\sigma^{o,y}$ represent the longitudinal and lateral standard deviations of the potential field distribution respectively. 

\subsection{Domain knowledge-based exploration branches}
\begin{figure} [h!]
    \centering
    \includegraphics[trim=12mm 0mm 5mm 0mm, clip, width=1.0\columnwidth]{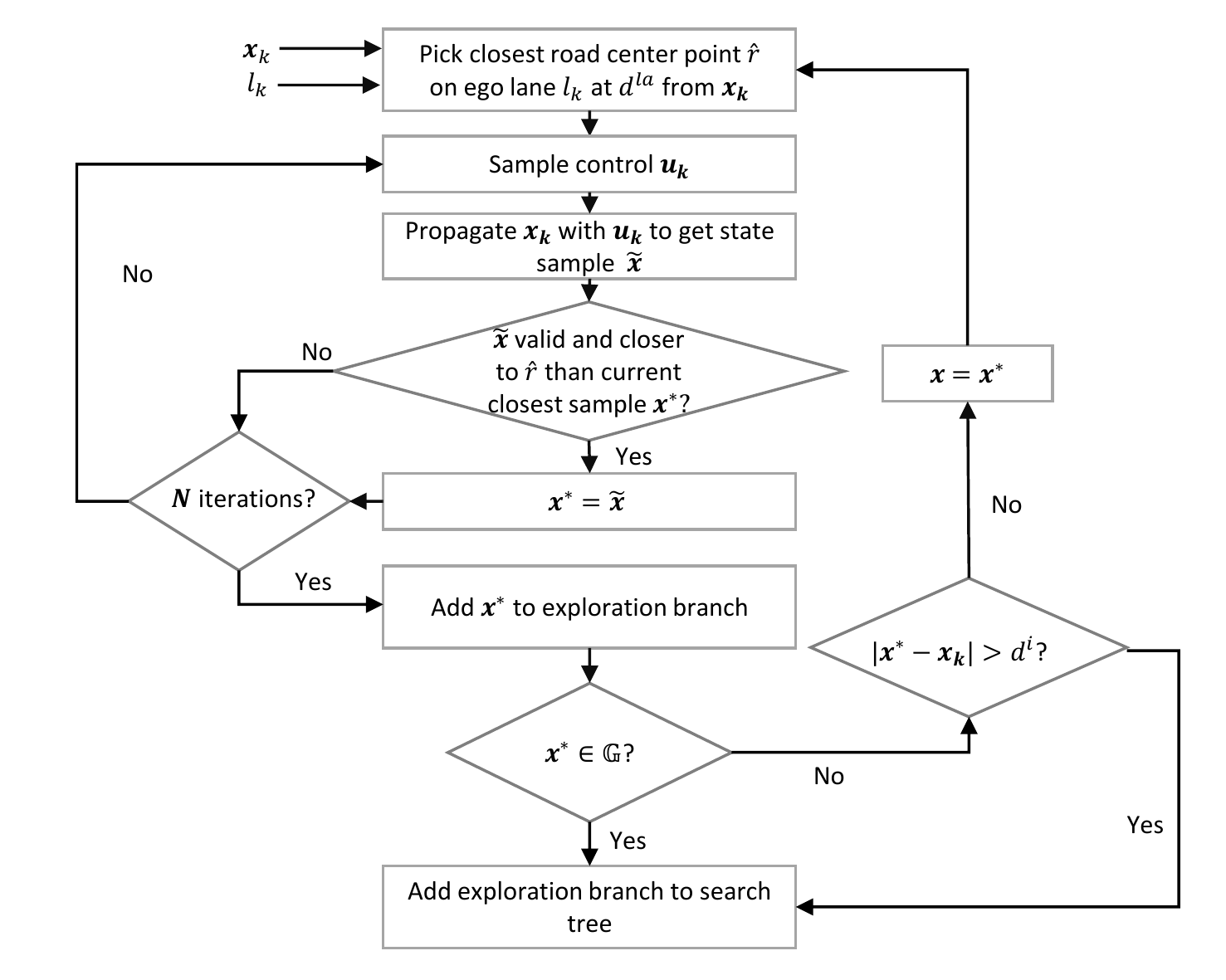}
    \caption{\footnotesize{Algorithm flowchart for adding road-layout based exploration branch.}}
    \label{fig:InitialSolution}
\end{figure}

The SST algorithm focuses on a trade-off between fast exploration of the available state space and fast convergence to high-quality paths. However, due to the lack of domain-specific knowledge on the state space, exploration can be slow or come at the expense of the final solution optimality. To overcome this, an informed exploration of the state space based on the road layout can be used. As the vehicle is expected to predominantly follow a (virtual) lane center while driving in a structured environment, an initial exploration branch of states is computed in the search tree, starting at tree-root $\mathbf{x}_0$, which steers the vehicle towards the lane center. This algorithm, for which a schematic overview is provided in Figure \ref{fig:InitialSolution}, picks a point on the current ego lane center at a fixed distance, $d^{la}$. Then, for a fixed number of $N$ iterations, inputs are sampled and applied to the current state to achieve a new state sample. After $N$ iterations, the sample closest to the lane center point, and corresponding control inputs, are saved. The initial state is set to the newly achieved state and the process is repeated until the goal space is reached, no valid state can be found or the distance from the starting state to the latest state sample exceeds a threshold.

\begin{figure}[h!]
    \centering
      \includegraphics[trim=10mm 0mm 7mm 0mm, clip, width=1.0\columnwidth]{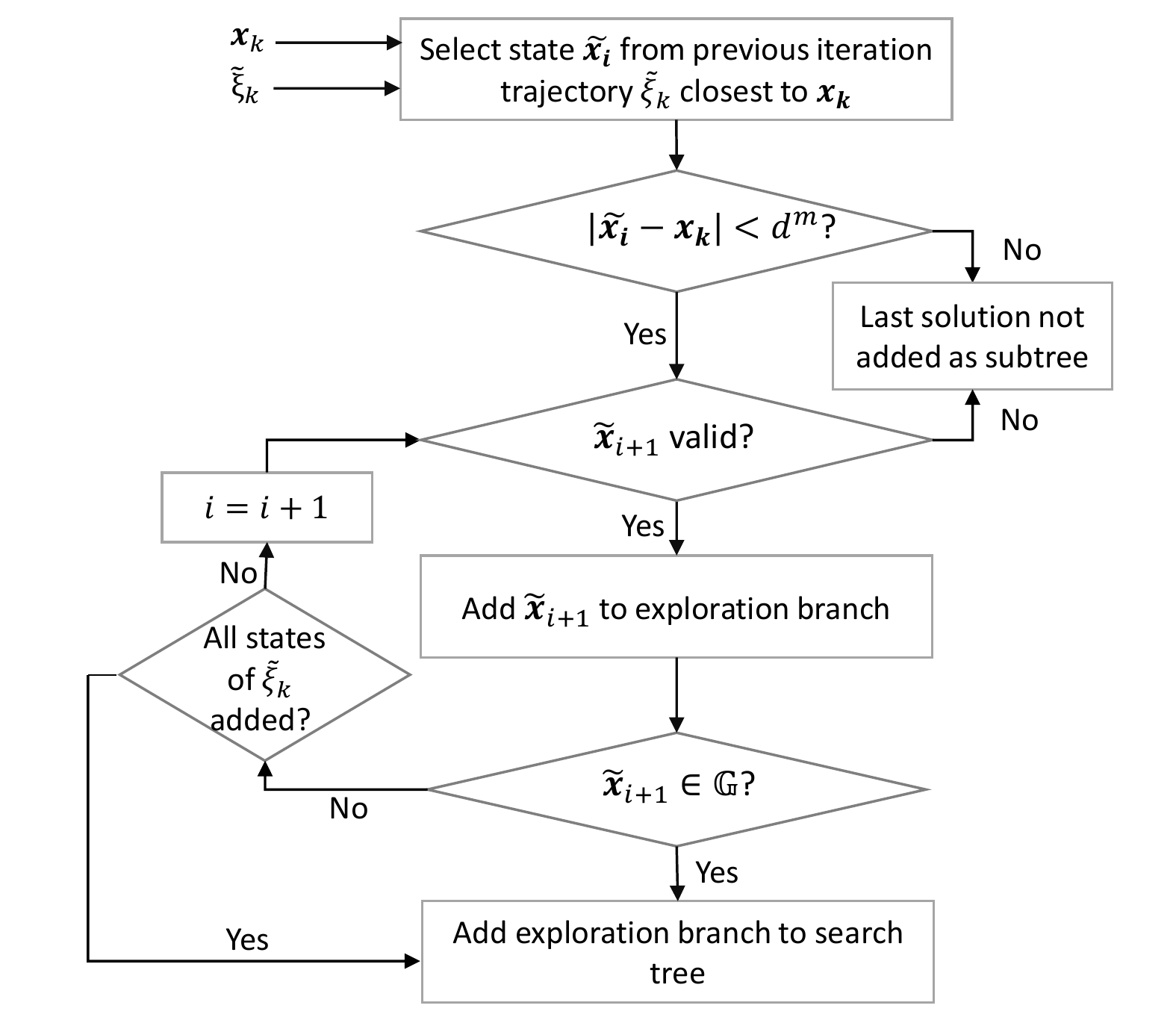}
    \caption{\footnotesize{Algorithm flowchart for adding the previous solution based exploration branch. The trajectory found in the previous iteration is denoted by $\tilde{\xi}$}}
    \label{fig:LastSolution}
\end{figure}

Furthermore, the original SST algorithm attempts to solve each planning query without knowledge of the previous planning iterations, i.e., previously computed trajectories and/or decisions are not taken into account. However, in case of a sufficiently high planner update rate and vehicle trajectory tracking capability, the vehicle is likely to be somewhere along the previously calculated trajectory. In this work, the search tree is augmented with a branch created from the trajectory calculated in the previous iteration, provided that the starting state, tree root $\mathbf{x}_0$, is sufficiently close to a given state on the (interpolated) previous solution. The schematic overview of the algorithm for embedding the computed trajectory of the previous planning iteration is depicted in Figure \ref{fig:LastSolution}. Here, the state of the previous trajectory is selected closest to the tree root. If this state is sufficiently close to the tree root, the states of the previous trajectory after the selected state are iteratively checked for validity and, if valid, added to the branch. If a state is invalid, the process is aborted. This does not only allow the planner to more quickly converge to high-quality paths as it can refine the previous iteration result, but it also raises the probability that any decision taken in a previous iteration (e.g., overtake versus the following predecessor) is adhered to in later iterations.S

The augmented planner as introduced in this section will henceforth be referred to as \textit{Domain-Knowledge-Informed-SST}, or in short \textit{DKI-SST}. The original algorithm will be referred to as \textit{base SST}.


\section{Scenario definition and validation}\label{sec:validation}
\begin{figure}[h!]
    \centering
    \includegraphics[trim=5mm 40mm 5mm 50mm, clip, width=1.0\columnwidth]{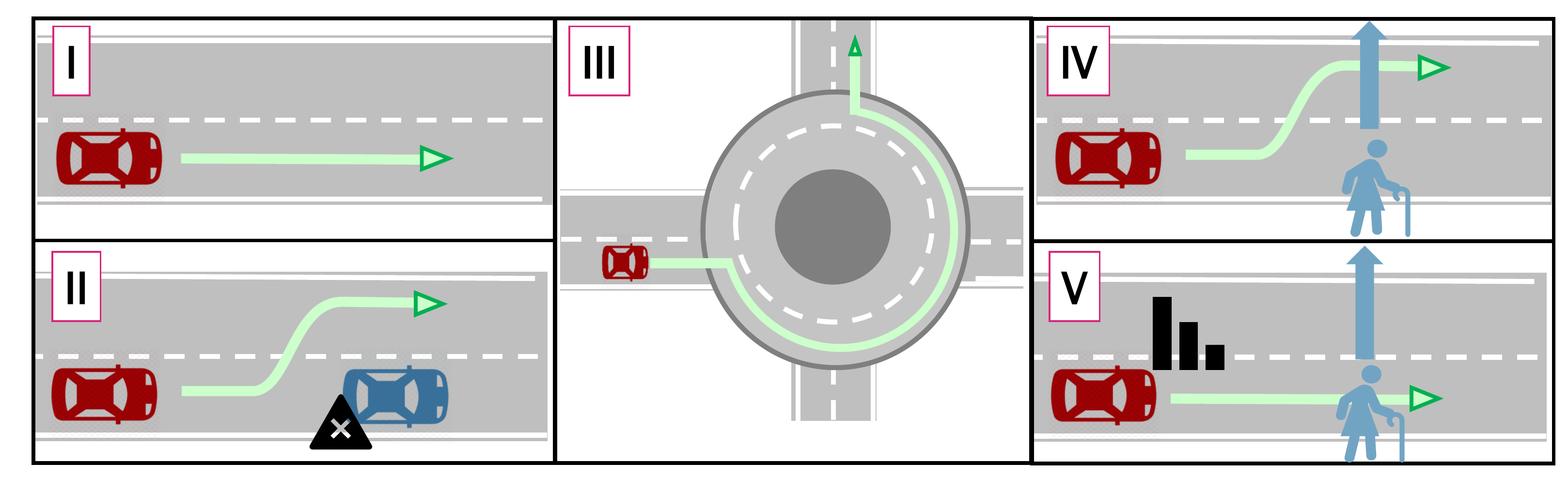}
    \caption{\footnotesize{Scenario definition for validation. I): following a straight road. II): overtaking a static predecessor. III): navigating a roundabout. IV): evasive maneuvering based on VRU motion prediction. V): speed adjustment based on VRU motion prediction}}
    \label{fig:scenarios}
\end{figure}

In this section, the method described in Section \ref{sec:method} will be proven based on several simulation studies. First, the performance of the DKI-SST planner is compared to the base SST planner in two simple scenarios to demonstrate the performance gain by introducing domain knowledge. Second, the performance of the DKI-SST planner is benchmarked based on three common urban driving scenarios, which include various road-layouts and dynamic traffic participants. A schematic overview of the scenarios is given in Figure \ref{fig:scenarios}. 

The automated vehicle is modeled as kinodynamic system with three degrees of freedom, as described in \cite{Pacejka2005}. This system is modeled using the Gazebo simulator \cite{Koenig2004}, where the road networks are described in the OpenDrive format \cite{Dupuis2010} and the dynamic actors are defined using the OpenScenario standard \cite{OpenScenario}. For vehicle actuation based on the computed trajectories, a feedback tracking controller is employed based on the deviation from the path and accompanying speed. The base SST planner is implemented in the OMPL framework \cite{OMPL2012}. 
The parametrization of the DKI-SST planner is depicted in Table \ref{table:pars}.

\begin{table}
\begin{tabular}{p{0.9cm} p{4.2cm} p{1.3cm} p{0.5cm}}\toprule
\textbf{Parameter}              & \textbf{Description} & \textbf{Value} & \textbf{Unit}\\
\toprule
$fq$ & Update rate & 2 & Hz\\
$t^q$ & Query solving time & 0.3 & s\\
$g^d$ & Fixed goal distance & 30.0 & m\\
$g^t$ & Longitudinal goal space threshold & 2.0 & m\\
$t^s$ & Integration time & 0.04 & s \\
$d^n$ & Nearest neighbour selection distance & 0.2 & - \\
$t^p$ & Propagation step & 0.4 & s\\
$d^p$ & Pruning distance & 0.1 & -\\
$\sigma^a$ & Acceleration input standard deviation & 0.8 & $\text{ms}^{-2}$\\
$\sigma^{\delta}$ & Steering angle input standard deviation & 0.2 & rad\\
$\underline{v},\overline{v}$&Min/max velocity&$[0.0,6.0]$&$\text{ms}^{-1}$\\
$\underline{a},\overline{a}$&Min/max acceleration&$[-0.8,0.8]$&$\text{ms}^{-2}$\\
$\underline{\delta},\overline{\delta}$&Min/max steering angle&$[-0.4,0.4]$&$\text{rad}$\\
$\overline{P}$ & Maximum penalty grid cell value & 100 & - \\
$\overline{p}$ & State invalid cell value & 99 & - \\
$w^{pl}$ & Path length penalty weight & 0.05 & - \\
$w^{ds}$ & Desired speed penalty weight & 0.5 & - \\
$w^{pg}$ & Penalty grid penalty weight & 0.2 & - \\
$w^{tc}$ & Target clearance penalty weight & 2.0 & - \\
$v^d$&Desired velocity&5.0&$\text{ms}^{-1}$\\
$A_O $ &Object risk amplitude& 100 &-\\
$\sigma^{o,x}$&Longitudinal risk deviation&3.0&$\text{m}$\\
$\sigma^{o,y}$&Lateral risk deviation&2.0&$\text{m}$\\
$d^{la}$&Road layout look ahead distance&3.0&$\text{m}$\\
$d^{i}$&Threshold distance for road layout exploration branch& 40.0 &$\text{m}$\\
$N$ & Max iterations for sampling road layout exploration state & 100 & - \\
$d^{m}$&Threshold distance between start state and previous iteration closest state & 1.0 & - \\
\toprule
\end{tabular}
\caption{\footnotesize{Simulation parameters for all considered cases.}}
\label{table:pars}
\end{table}

\subsection{Comparison to base SST planner}
The performance gained by introducing domain knowledge is demonstrated based on two simple scenarios: driving on a straight road and overtaking a static object. For the first case, the vehicle should stay in the middle of the lane as much as possible while both longitudinal and lateral accelerations should be limited to drive comfortably. Furthermore, the vehicle should drive as close as possible to the desired velocity, as it is not limited by safety-critical environmental factors. To circumvent the possible inaccuracies of a tracking controller in this scenario, the vehicle is manually driven at exactly the lane center, at exactly the desired speed ($5m/s$). For each calculated trajectory, the mean absolute values for the deviation from the desired velocity, accelerations and deviation from closest lane center are taken. The simulation is run for a period of minimum 10 sec. 
For the second  considered case, when overtaking a static object, it is important that the vehicle initiates a timely overtaking maneuver, thus preventing coming too close to the object or causing a collision. Therefore, the minimal distance to the object is taken as a performance indicator from this scenario. As vehicle actuation is required in this scenario, the tracking controller is tuned for optimal performance for the Base SST and DKI-SST algorithms separately. 

The overview of the results of both planners, shown in Table \ref{table:comparisonToBase} clearly show the benefits of the DKI-SST method. This leads to 75\% performance gain in the use of acceleration and more than 95\% performance gain in the deviation from desired speed. 
\begin{table}[h!]
\small
\begin{tabular}{ p{3.4cm} p{1.0cm} p{1.0cm} p{1.5cm}}
\toprule
\textbf{Performance criteria} & \textbf{Base SST} & \textbf{DKI-SST} & \textbf{Performance gain [\%]} \\
\toprule
 Mean absolute trajectory longitudinal acceleration $[\text{ms}^{-2}]$ & 0.52 & 0.13 & 75.0  \\ 
 Mean absolute trajectory deviation from desired speed $[\text{ms}^{-1}]$ & 1.52 & 0.051 & 96.6 \\
 Mean absolute trajectory deviation from closest lane center $[\text{m}]$ & 0.33 & 0.057 & 82.7 \\
 Minimum resulting distance to target $[\text{m}]$ & 3.87 & 3.94 & 1.8\\
 \toprule
\end{tabular}
\caption{\footnotesize{Performance indicators for the base SST algorithm versus DKI-SST.}}
\label{table:comparisonToBase}
\end{table}
\subsection{Roundabout}
\begin{figure}[h!]
    \centering
    \includegraphics[trim=0mm 0mm 0mm 0mm, clip, scale=0.25]{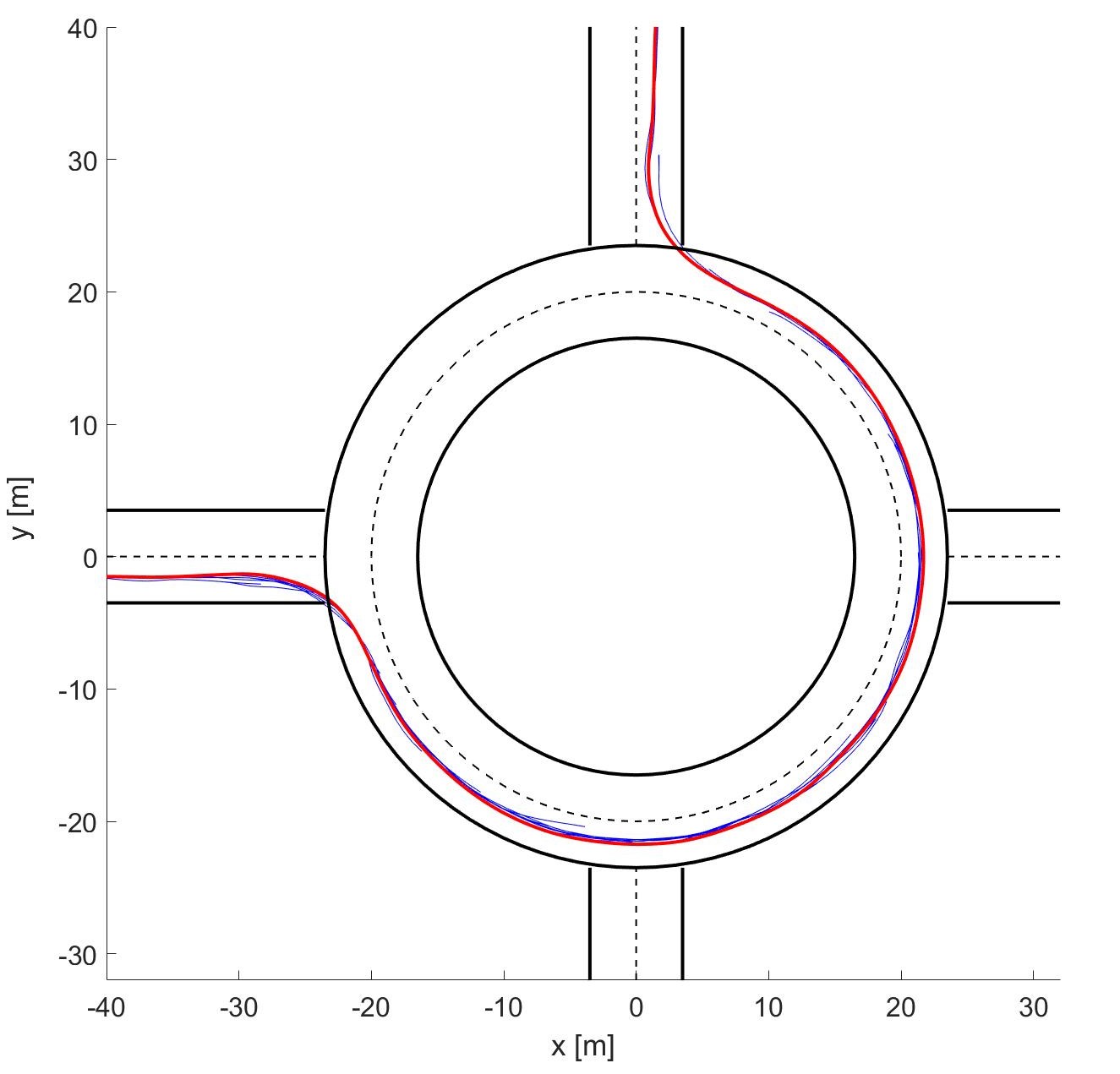}
    \caption{\footnotesize{Planned trajectories (blue) and final driven path (red) on a roundabout}}
    \label{fig:UC3TrajectoriesAlt}
\end{figure}

The third scenario under consideration is navigating a roundabout. The discrete change in heading angle of the road when entering the roundabout, together with the smooth curve on the roundabout itself, make this a challenging road layout for urban driving. In this scenario, the vehicle is guided to take the roundabout at the third exit. In Figure  \ref{fig:UC3TrajectoriesAlt}, it is shown that the ground plane positions of all computed trajectories (blue) are navigating the vehicle within the lane over the roundabout. This results in a safe and comfortable driven path, depicted in red. 
\subsection{Preventive steering for Vulnerable Road Users}
\begin{figure}[h!]
        \centering
        \begin{subfigure}[b]{0.90\columnwidth}
            \centering
            \includegraphics[trim=30mm 0mm 30mm 0mm, clip, width=\columnwidth]{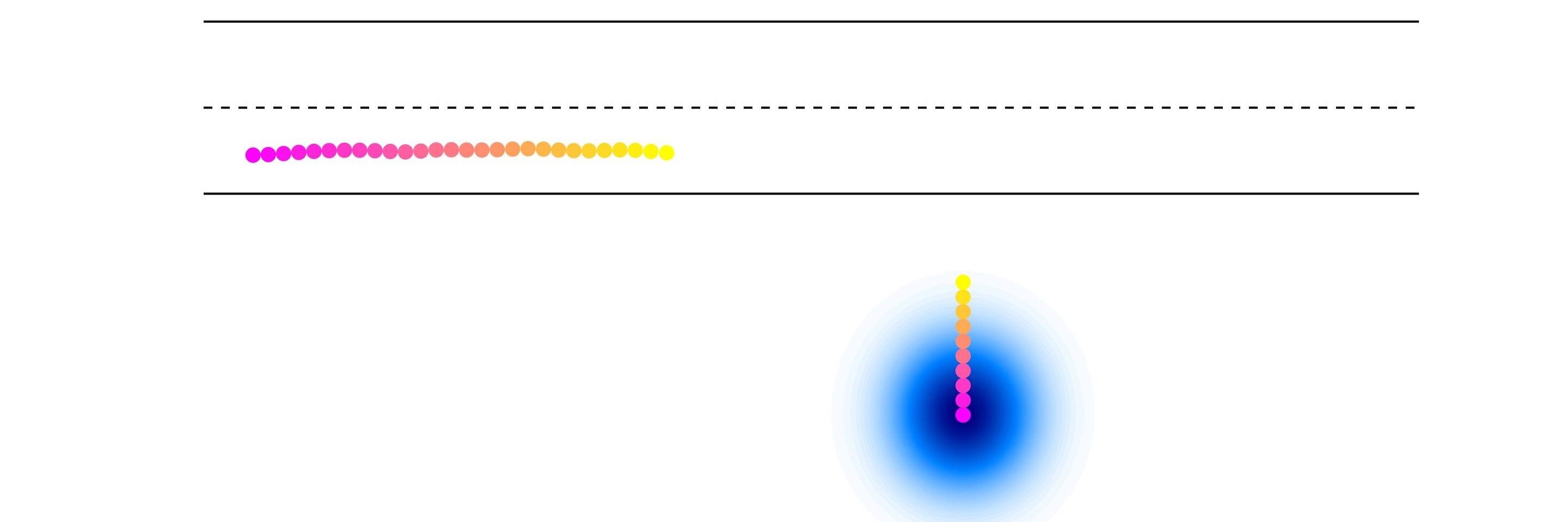}
            \caption[]%
            {{\footnotesize t=51.8s}}    
            \label{fig:27b_1}
        \end{subfigure}
        \begin{subfigure}[b]{0.95\columnwidth}  
            \centering 
            \includegraphics[trim=30mm 0mm 30mm 0mm, clip, width=\columnwidth]{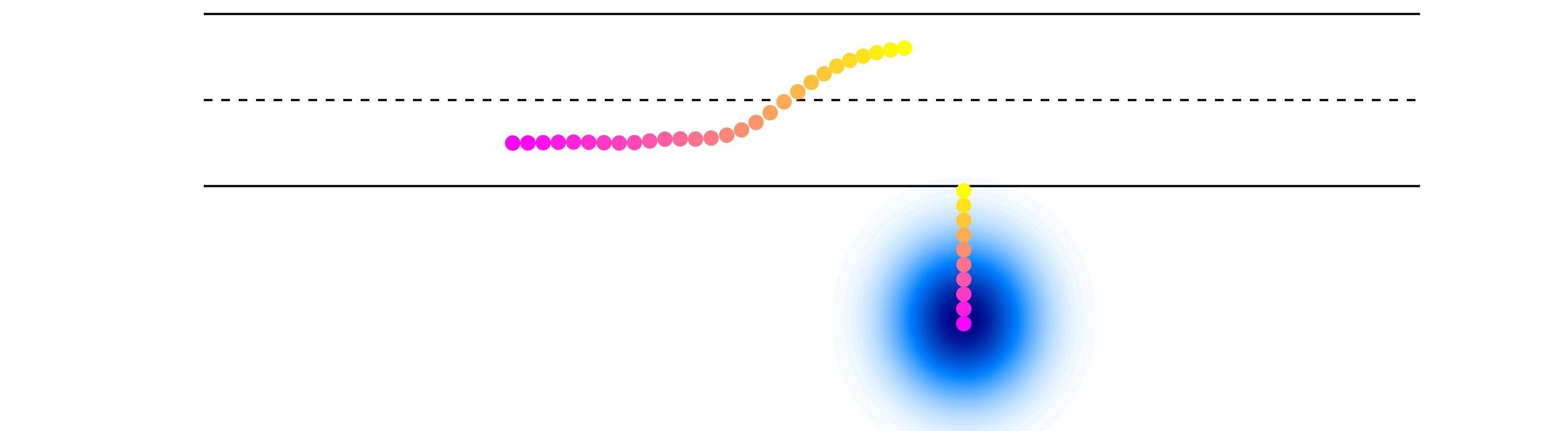}
            \caption[]%
            {{\footnotesize t=55.3s}}    
            \label{fig:27b_2}
        \end{subfigure}
        \vskip\baselineskip
        \begin{subfigure}[b]{0.95\columnwidth}   
            \centering 
            \includegraphics[trim=30mm 0mm 30mm 0mm, clip, width=\columnwidth]{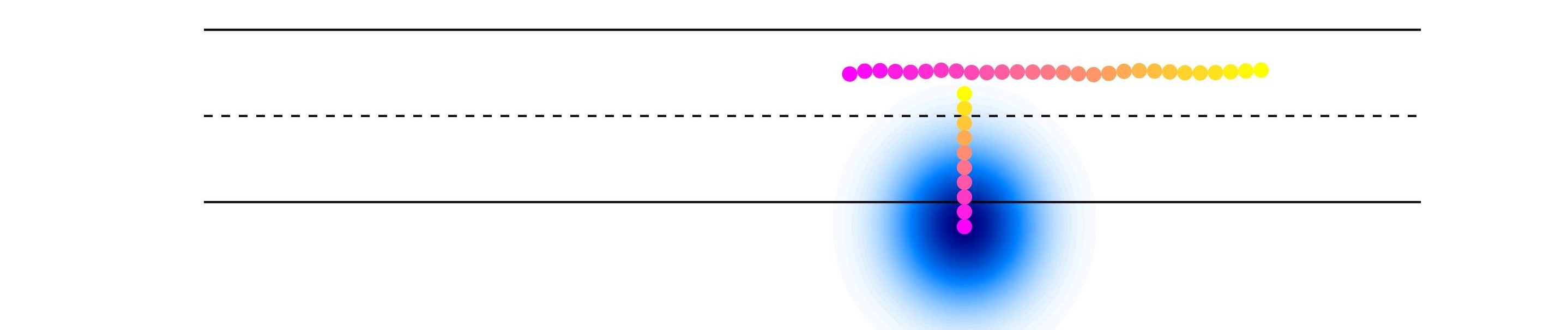}
            \caption[]%
            {{\footnotesize t=59.8s}}    
            \label{fig:27b_3}
        \end{subfigure}
        \hfill
        \begin{subfigure}[b]{0.95\columnwidth}   
            \centering 
            \includegraphics[trim=30mm 0mm 30mm 0mm, clip, width=\columnwidth]{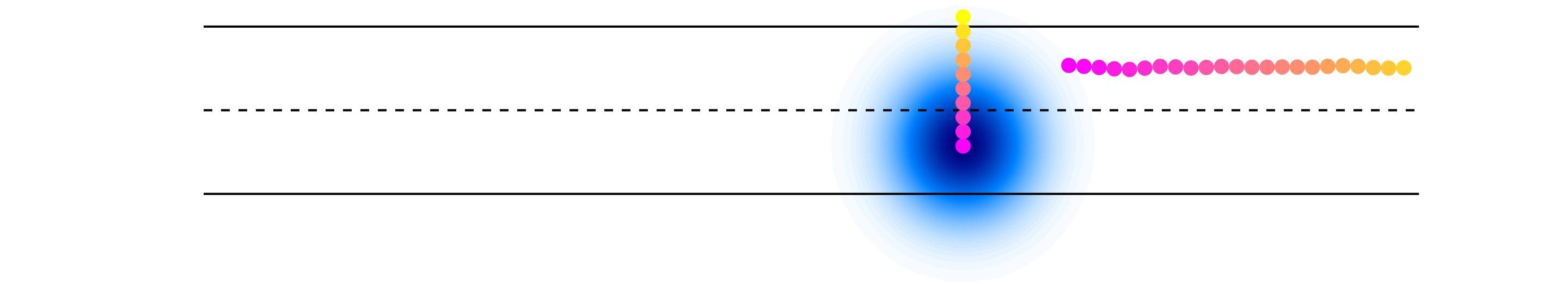}
            \caption[]%
            {{\footnotesize t=62.8s}}    
            \label{fig:27b_d}
        \end{subfigure}
        \caption[]
        {\footnotesize{Planned trajectory of the vehicle and predicted trajectory of the pedestrian, including its current repulsive potential field (blue). Trajectories range from current time (purple) to 5 seconds ahead (yellow)}. }
        \label{fig:27b}
\end{figure}

To determine the trajectory planners capability of accounting for object motion prediction, two VRU-crossing scenarios are evaluated. In scenario \rom{4}, the VRU has merely started crossing the two-lane street when the automated vehicle approaches. In this case, a preemptive lane change is a safe and valid option. The resulting calculated trajectories are depicted in Figure \ref{fig:27b} at four different time instances. Here, it is shown that the vehicle makes are comfortable lane change, well before the current VRU repulsive field requires to do so. The preemptive lane change is initiated due to the expected trajectory of the VRU, which is expected to be at a sufficient lateral distance for the vehicle to safely pass after a lane change. 

\subsection{Preventive braking for Vulnerable Road Users}
\begin{figure}
        \centering
        \begin{subfigure}[b]{0.95\columnwidth}
            \centering
            \includegraphics[trim=30mm 0mm 30mm 0mm, clip, width=\columnwidth]{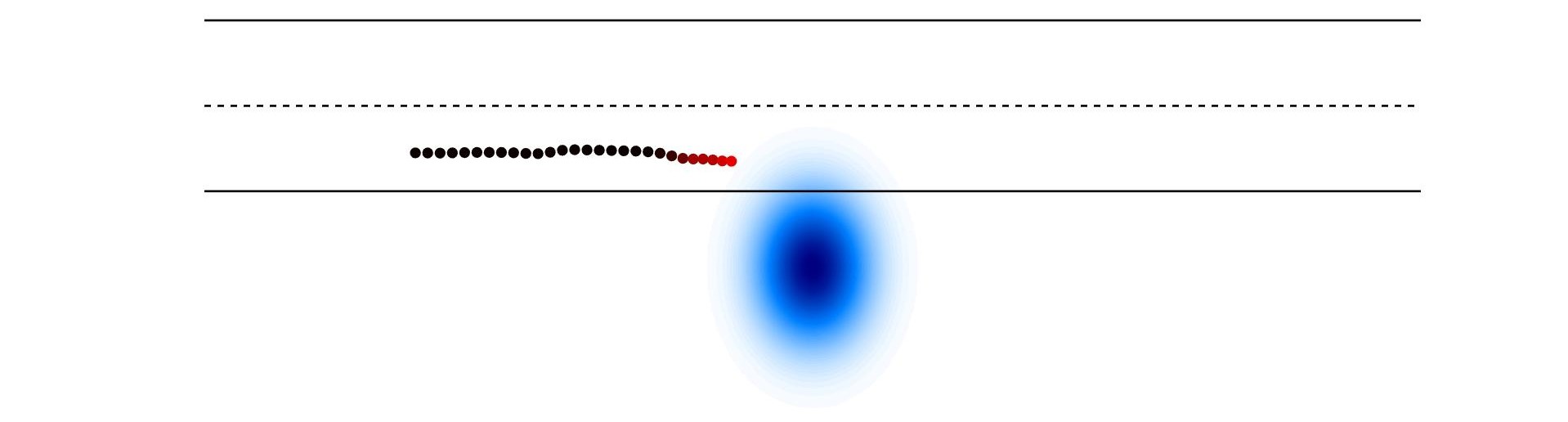}
            \caption[]%
            {{\footnotesize t=49.5s}}    
            \label{fig:27c_1}
        \end{subfigure}
        \begin{subfigure}[b]{0.95\columnwidth}  
            \centering 
            \includegraphics[trim=30mm 0mm 30mm 0mm, clip, width=\columnwidth]{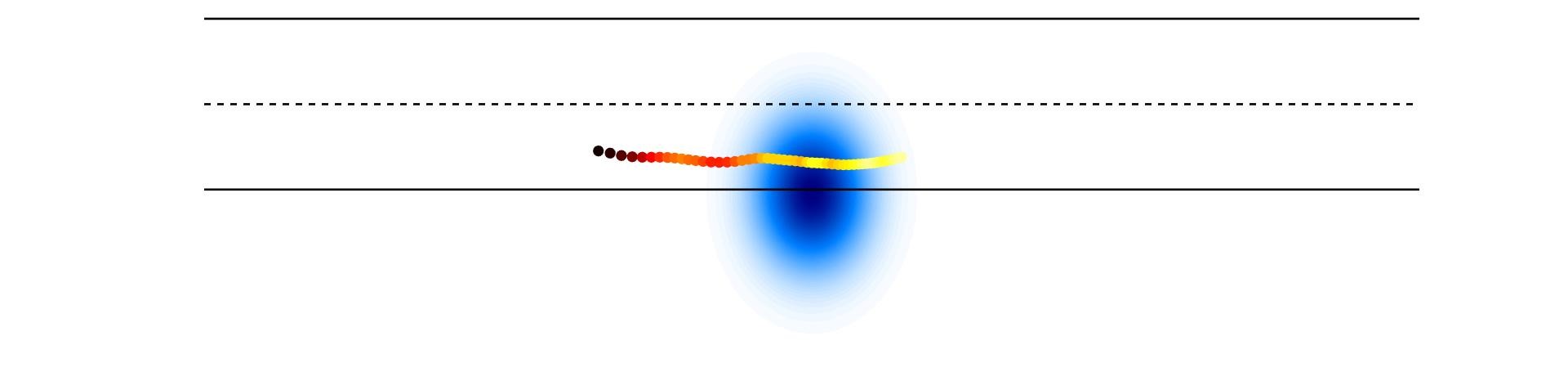}
            \caption[]%
            {{\footnotesize t=52.5s}}    
            \label{fig:27c_2}
        \end{subfigure}
        \vskip\baselineskip
        \begin{subfigure}[b]{0.95\columnwidth}   
            \centering 
            \includegraphics[trim=30mm 0mm 30mm 0mm, clip, width=\columnwidth]{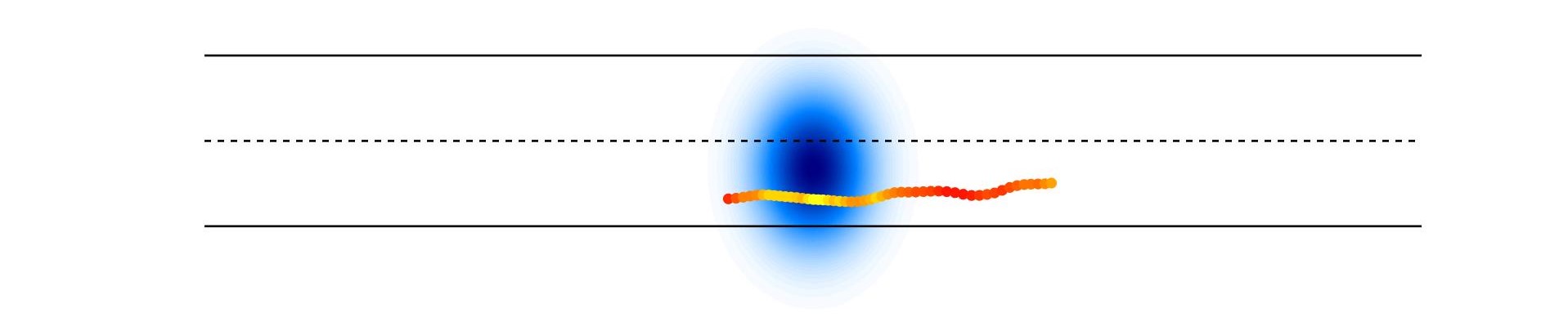}
            \caption[]%
            {{\footnotesize t=55.0s}}    
            \label{fig:27c_3}
        \end{subfigure}
        \hfill
        \begin{subfigure}[b]{0.95\columnwidth}   
            \centering 
            \includegraphics[trim=30mm 0mm 30mm 0mm, clip, width=\columnwidth]{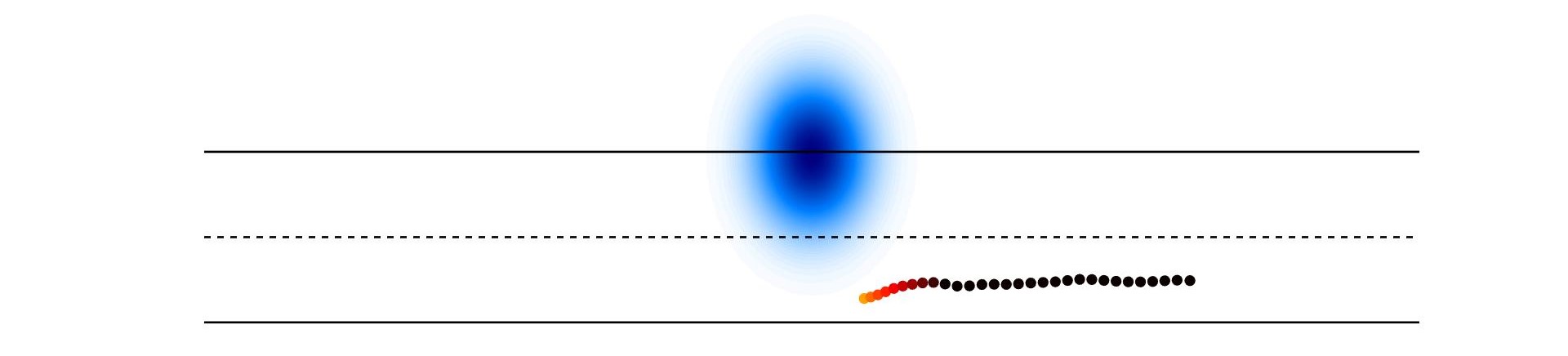}
            \caption[]%
            {{\footnotesize t=59.5s}}    
            \label{fig:27c_d}
        \end{subfigure}
        \caption[]
        {\footnotesize{Planned trajectory of the vehicle, together with the current repulsive potential field of the pedestrian. The coloring of the trajectory indicates the computed speed at that point, ranging from black (5$m/s$) to yellow (1$m/s$)}.  }
        \label{fig:27c}
\end{figure}

\begin{figure}[h!]
    \centering
    \includegraphics[width=0.9\columnwidth]{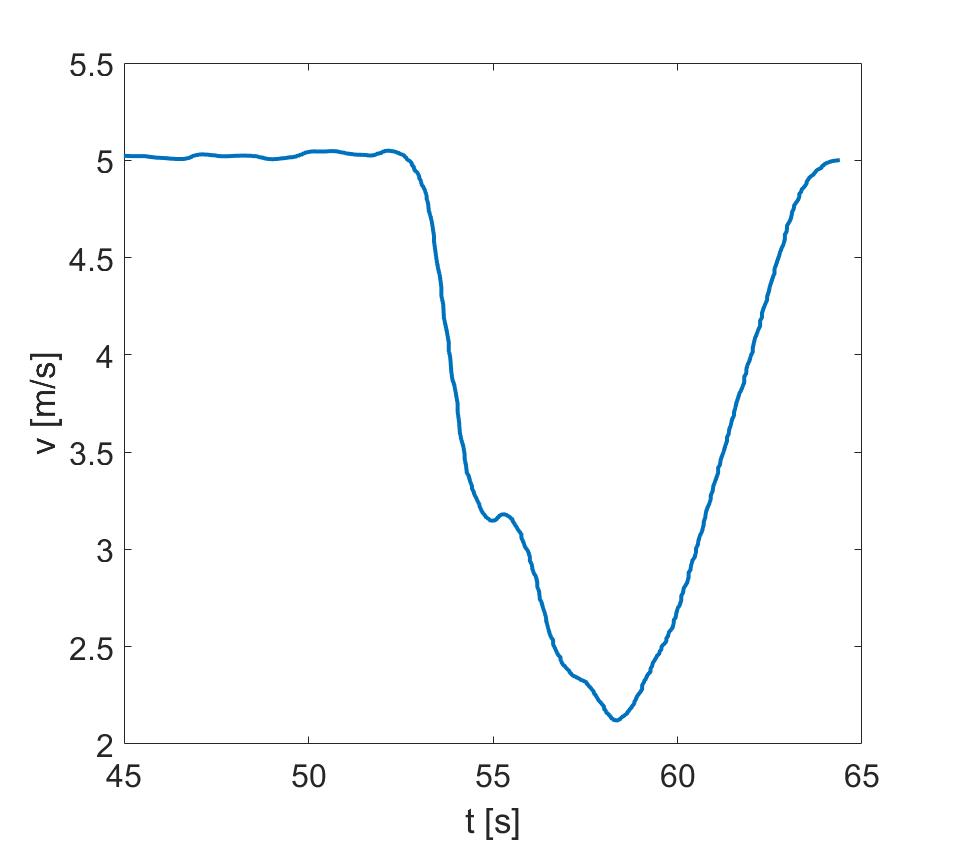}
    \caption{\footnotesize{Longitudinal vehicle velocity over time of scenario \rom{5}}}
    \label{fig:UC27cVelocity}
\end{figure}

In the second VRU-crossing scenario, the VRU is well into crossing the road when the vehicle approaches. Due to the movement of the VRU, solely evasive steering is not a safe option. Figure \ref{fig:27c} shows the computed paths again at four different time instances, together with the velocity profile along this path. The vehicle reduces its speed along the computed trajectories, to allow the VRU to safely cross the road, before speeding up again. The resulting vehicle speed is shown in Figure \ref{fig:UC27cVelocity}. 

\section{Conclusions}\label{sec:conclusion}
In this work, a real-time and high-quality trajectory planner is introduced for safe and comfortable automated driving in dynamic urban environments which takes explicitly into account the road layout and the predicted motion of other traffic participants. Building on a previous SST algorithm, additional road-layout and previous trajectories knowledge are embedded to optimize the generation and exploration of the search tree branches.
Results show significant performance gains of this new method when compared to the baseline on several criteria. 
Furthermore, using several urban scenarios, the planner has shown to be able to navigate various road layouts, as well as be able to adapt both planned path and speed profile to the predicted motion of other traffic participants. 
Future work includes shaping the repulsive artificial potential function based on object information to include even more realistic information, e.g. estimation uncertainty and object type, to include more complex context awareness and validate these through real-life experiments.  

\section{Acknowledgments}\label{sec:acknowledgements}
This work is supported by the EU Horizon 2020 R{\&}D program under grant agreement No. 861570, project SAFE-UP (proactive SAFEty systems and tools for a constantly UPgrading road environment). 

\printbibliography

\end{document}